\documentclass[10pt]{article}
\newfont{\Bb}{msbm10 scaled\magstep1}
\newfont{\Bbs}{msbm10 scaled\magstep0}

\newtheorem{theorem}{Theorem}[section]

\newtheorem{lemma}[theorem]{Lemma}

\newtheorem{proposition}[theorem]{Proposition}
\newtheorem{definition}[theorem]{Definition}
\newtheorem{conjecture}[theorem]{Conjecture}

\newenvironment{exafont}{\begin{bf}}{\end{bf}}
\newenvironment{example}{\vspace{0.3cm}\par\noindent\refstepcounter{theorem}\begin{exafont}Example \thetheorem\end{exafont}\hspace{\labelsep}}{\vspace{0.3cm}\par}

\def\F{{\bf F}_q}

\title{Rationality of Partial Zeta Functions}

\author{Daqing Wan\footnote{Partially supported by the NSF and the NSFC}}

\begin{document}

\maketitle

\begin{abstract}
We prove that the partial zeta 
function introduced in \cite{W3} is a rational function, 
generalizing Dwork's rationality theorem. 
\end{abstract}

\section{Introduction}

Let $\F$ be the finite field of $q$ elements of characteristic $p$.   
Let $\bar{\bf F}_q$ be a fixed algebraic closure of $\F$. 
Let $X$ be an affine algebraic variety over $\F$, 
embedded in some affine space ${\bf A}^n$.   
That is, $X$ is defined by a system of polynomial equations
$$F_1(x_1,\cdots, x_n)=\cdots =F_m(x_1,\cdots, x_n)=0,$$
where each $F_i$ is a polynomial 
defined over $\F$. Let $d_1,\cdots, d_n$ be positive integers. 
For each positive integer $k$, let 
$$X_{d_1,\cdots, d_n}(k)=\{ x\in X(\bar{\bf F}_q)| 
x_1\in {\bf F}_{q^{d_1k}}, \cdots, x_n\in {\bf F}_{q^{d_nk}}\}.$$ 
The number $\# X_{d_1,\cdots, d_n}(k)$ counts the points of $X$ 
whose coordinates are in different subfields of $\bar{\bf F}_q$. 
We would like to understand this sequence of 
integers $\# X_{d_1,\cdots, d_n}(k)$ indexed by $k$.  
As usual, it is sufficient to understand the following generating 
function.

\begin{definition} Given $X$ and the $n$ positive integers $d_1,\cdots, d_n$, 
the associated partial zeta function $Z_{d_1,\cdots, d_n}(X,T)$  
of $X$ is defined to be the following formal power series
$$Z_{d_1,\cdots, d_n}(X, T) = \exp \left( \sum_{k=1}^{\infty} 
{\# X_{d_1,\cdots, d_n}(k) \over k}T^k \right) \in 1+T{\bf Q}[[T]].$$
\end{definition}
 
Replacing $q$ by a power of $q$, without loss of generality, we may assume that 
the integer $d_i$'s are relatively prime. 
In the special case that $d_1=\cdots =d_n=1$, the number $\# X_{1, \cdots, 1}(k)$ 
is just the number of ${\bf F}_{q^{k}}$-rational points on $X$. The 
partial zeta function $Z_{1,\cdots, 1}(X,T)$ then becomes the 
classical zeta function $Z(X, T)$ of the variety 
$X$.  Dwork's rationality theorem \cite{Dw} says that $Z(X,T)$ 
is a rational function. Deligne's theorem \cite{De} on the Weil conjectures 
says that the reciprocal zeros and the reciprocal poles of $Z(X, T)$ 
are Weil $q$-integers. Recall that a Weil $q$-integer $\alpha$ is an 
algebraic integer such that $\alpha$ and each of its Galois conjugates 
have the same complex absolute value $q^{w/2}$ for some non-negative integer $w$. 
The integer $w$ is called the weight of $\alpha$. 

One of our motivations to introduce the above more general partial zeta function 
comes from potential applications in number theory, combinatorics and 
coding theory. From a theoretic point of view, a special case of the partial 
zeta function reduces to the geometric moment zeta function \cite{W4} attached to a 
family of algebraic varieties over $\F$, which was in turn motivated by our work on   
Dwork's unit root conjecture \cite{W1}\cite{W2}. Intuitively, the partial zeta function 
gives many new ways, parametrized by the integers $d_i$'s, to count the geometric 
points on $X$ and thus it contains critical information about the distribution 
of the geometric points of $X$. The partial zeta function also provides 
a simple diophantine reformulation 
of many much more technical problems. In \cite{W3}, the following two results were 
proven concerning the possible rationality of the partial zeta function. 
Recall that for a complex number $\alpha$ and 
a complex power series $R(T)$ with constant term $1$, we can define 
the complex power $R(T)^{\alpha}$ in terms of 
the binomial series $(1 + T{R(T)-1\over T})^{\alpha}$. 

\begin{proposition}[Faltings \cite{W3}] Let $d=[d_1,\cdots, d_n]$ be the least 
common multiple of the $d_i$. Let $\zeta_d$ be a primitive $d$-th root 
of unity.  There are $d$ rational functions 
$R_j(T)$ ($1\leq j\leq d$) with $R_j(0)=1$ and with algebraic 
integer coefficients such that 
$$Z_{d_1,\cdots, d_n}(X,T)=\prod_{j=1}^d R_j(T)^{\zeta_d^j}.$$
Furthermore, the reciprocal zeros and reciprocal poles of 
the $R_j(T)$'s are Weil $q$-integers. 
\end{proposition}

This result shows that the partial zeta function 
is nearly rational. It is proved by using a geometric construction of Faltings and 
the general fixed point theorem in $\ell$-adic cohomology.

\begin{proposition}[\cite{W3}] If the integers $\{ d_1, d_2,\cdots, d_n\}$ 
can be rearranged such that $d_1|d_2|\cdots |d_n$, then 
the partial zeta function $Z_{d_1,\cdots, d_n}(X,T)$ is a rational 
function in $T$, who reciprocal zeros and reciprocal poles are Weil $q$-integers.  
\end{proposition}

This result shows that the partial zeta function has the stronger property 
of being a rational function in some non-trivial special cases. 
It is proved by viewing $X$ 
as a sequence of fibered varieties and inductively using the Adams operation of the relative 
$\ell$-adic cohomology. Although it was felt that the partial zeta function 
may not be always rational in general, no counter-examples were found.  
The aim of this note is to prove the following result. 
 
\begin{theorem} For any variety $X$ as above and any positive integers $\{ d_1,\cdots, d_n\}$,  
the partial zeta function $Z_{d_1,\cdots, d_n}(X,T)$ is 
a rational function in $T$,  who reciprocal zeros and reciprocal poles are 
Weil $q$-integers.
\end{theorem}

The idea of the proof is to exploit the geometric construction of Faltings 
and its relation to Galois action.   
Once the rationality is proved, 
one main new problem about the partial zeta function is to understand 
its dependence and variation on the arithmetic 
parameters $d_i$'s. This would raise many interesting new questions to be explored,  
as already illustrated in the special case of moment zeta functions \cite{W4}.  
The first question one could ask for is about the number of zeros and poles of the partial 
zeta function. In Fu-Wan \cite{FW1}, using Katz's bound \cite{Ka1} on the $\ell$-adic Betti numbers, 
an explicit total degree bound for $Z_{d_1,\cdots, d_n}(X,T)$ is given, which grows exponentially 
in $d$. We conjecture that the true size of the total degree is much smaller and 
bounded by a polynomial in $d$. 

\begin{conjecture} There are two positive constants $c_1(X)$ and $c_2(X)$ depending only 
on $X$  such that the total degree of the partial zeta function 
$Z_{d_1,\cdots, d_n}(X,T)$ is uniformly 
bounded by $c_1(X)d^{c_2(X)}$ for all positive integers $\{ d_1,\cdots, d_n\}$.  
\end{conjecture}

This conjecture has been proved to be true in Fu-Wan \cite{FW2} in the special case that 
$d_1=\cdots = d_r =1$ and $d_{r+1}=\cdots =d_n =d$, corresponding to the so-called 
moment zeta function case which has been studied more extensively in connection 
to Dwork's unit root conjecture. We believe that the above conjecture (if true) 
together with a deeper analysis of the weights of the zeros and poles of the partial 
zeta function would have many important applications. Under suitable conditions, 
we would like to have optimal estimates of the form 
$$|\# X_{d_1,\cdots, d_n}(k) -q^{k(d_1+\cdots +d_n-dm)}| \leq c_1d^{c_2} 
q^{k(d_1+\cdots +d_n-dm)/2},$$
see section $4$ for some results in the case of 
Artin-Schreier hypersurfaces ($m=1$).

{\bf Acknowledgements}. 
Some results of this paper were obtained during 
the 2001 Lorentz center workshop ``L-functions from algebraic geometry" at Leiden University  
and the 2003 AIM workshop ``Future directions in algorithmic number theory". 
The author thanks both institutes for their hospitality. 
The author would also like to thank H. W. Lenstra Jr. for his interests 
and discussions on this paper.

\section{Rationality of partial zeta functions }

We slightly generalize the setup in the introduction. 
Let $f_i: X\to X_i$ ($1\leq i\leq n$) be morphisms of schemes of finite type over $\F$. 
Assume that the map $f: X \to X_1\times \cdots \times X_n$ defined 
by 
$$f(x) = (f_1(x), \cdots, f_n(x))$$ 
is an embedding. 
For each positive integer $k$, let 
$$f_{d_1,\cdots, d_n}(k)=\{ x\in X(\bar{\bf F}_q)| 
f_1(x)\in X_1({\bf F}_{q^{d_1k}}), \cdots, f_n(x)\in X_n({\bf F}_{q^{d_nk}})\}.$$  
This is a finite set since $f$ is an embedding. 

\begin{definition} Given the morphism $f$ and the $n$ positive integers $\{d_1,\cdots, d_n\}$, 
the associated partial zeta function $Z_{d_1,\cdots, d_n}(f,T)$ 
of the morphism $f$ is defined to be the following formal power series
$$Z_{d_1,\cdots, d_n}(f, T) = \exp \left( \sum_{k=1}^{\infty} 
{\# f_{d_1,\cdots, d_n}(k) \over k}T^k \right) \in 1+T{\bf Q}[[T]].$$
\end{definition}

It is clear that the special case in the introduction corresponds to the case that 
$X$ is affine in ${\bf A}^n$ and $f_i$ is the projection 
of $x$ to the $i$-th coordinate $x_i\in {\bf A}^1$. 

\begin{theorem} For any morphism $f$ and any positive integers $\{d_1,\cdots, d_n\}$, 
the partial zeta function $Z_{d_1,\cdots, d_n}(f,T)$ is  
a rational function in $T$, whose reciprocal zeros and reciprocal poles 
are Weil $q$-integers.  
\end{theorem}

To prove this theorem, we begin with the geometric construction of Faltings. 
Let $d=[d_1,\cdots, d_n]$ be the least common multiple. 
The set of geometric points on the $d$-fold product $X^d$ of $X$ 
has two commuting actions. One is the $q^{-1}$-th power 
geometric Frobenius action 
denoted by ${\rm Frob}$. Another is the automorphism $\sigma$ on $X^d$ 
defined by the cyclic shift 
$$\sigma(y_1,\cdots, y_d)=(y_d, y_1,\cdots, y_{d-1}),$$
where $y_j$ denotes the $j$-th component ($1\leq j \leq d$) 
of a point $y=(y_1,\cdots, y_d)$ on 
the $d$-fold product $X^d$. Thus, each component 
$y_j$ is a point on $X$.  
Let $Y=Y(d_1,\cdots, d_n,f)$ be the subvariety of $X^d$ 
defined by the equations   
$$f_i\circ \sigma^{d_i}=f_i, \ 1\leq i\leq n,$$
where $f_i: X^d \to X_i^d$ denotes the map $f_i(y_1,\cdots, y_d)=(f_i(y_1),\cdots, f_i(y_d))$. 
Thus, a point $y=(y_1,\cdots, y_d)\in X^d$ is on the 
subvariety $Y$ if and only if 
$$f_i(y_{j}) =f_i(y_{j+d_i}), \ 1\leq i\leq n, \ 1\leq j\leq d, \eqno(2.1)$$
where $j+d_i$ is taken to be the smallest positive residue of 
$j+d_i$ modulo $d$. It is clear that $Y$ is stable under the action of 
$\sigma$ which commutes with ${\rm Frob}$.

Now, let $a$ be a fixed positive integer relatively prime to $d$. 
Let $y=(y_1,\cdots, y_d)$ be a geometric point of $Y$. 
One checks that 
$$\sigma^a \circ {\rm Frob}^k(y)=y \Longleftrightarrow 
{\rm Frob}^k(y_{j})=y_{j+a}, \  1\leq j\leq d. \eqno(2.2)$$
The latter is true if and only if 
$${\rm Frob}^k(f_i(y_{j}))=f_i(y_{j+a}), \ 1\leq i\leq n, \ 1\leq j\leq d \eqno(2.3)$$
as $f$ is an embedding. 
Iterating equation (2.3) $d_i$ times, we get 
$${\rm Frob}^{d_ik}(f_i(y_{j}))=f_i(y_{j+ad_i}).$$
Since $y$ is on $Y$, by (2.1), we deduce that 
$${\rm Frob}^{d_ik}(f_i(y_{j}))=f_i(y_{j}).$$
Taking $j=1$, we see that 
every fixed point $y\in Y(\bar{\bf F}_q)$ of $\sigma^a \circ {\rm Frob}^k$ 
uniquely determines a point $y_1\in X(\bar{\bf F}_q)$ 
satisfying $f_i(y_1)\in X_i({\bf F}_{q^{d_ik}})$ 
for all $1\leq i\leq n$.    

Conversely, given $y_1\in X(\bar{\bf F}_q)$ such that $f_i(y_1)\in X_i({\bf F}_{q^{d_ik}})$ 
for all $1\leq i\leq n$, we define 
$$y_j = {\rm Frob}^{kh_j}(y_1), \ 1\leq j\leq d, $$
where $h_j$ is the unique integer between $0$ and $d-1$ such that 
$ah_j+1 \equiv j (~{\rm mod}~d)$. The integer $h_j$ is clearly well defined 
since $a$ and $d$ are relatively prime. 
If $j\equiv j' (~{\rm mod}~d_i)$, then $h_j\equiv h_{j'} (~{\rm mod}~d_i)$. 
Since $f_i(y_1)\in X_i({\bf F}_{q^{d_ik}})$, we deduce that 
$$f_i(y_j) = {\rm Frob}^{kh_j}(f_i(y_1)) = {\rm Frob}^{kh_{j'}}(f_i(y_1)) = f_i(y_{j'}).$$
This shows that the point $y=(y_1,\cdots, y_d)$ is on $Y$. 
Since $f$ is an embedding and $f_i(y_1)\in X_i({\bf F}_{q^{d_ik}})$ for all $i$, we deduce 
that $y_1 \in X({\bf F}_{q^{dk}})$. Using the congruence $a(h_j+1)+1 \equiv j+a (~{\rm mod}~d)$, 
we derive that 
$${\rm Frob}^k(y_j) ={\rm Frob}^{k(h_j+1)}(y_1) = y_{j+a}.$$
This proves that $\sigma^a \circ {\rm Frob}^k(y)=y$. 
In summary, we have proved the following result.

\begin{lemma} Let $a$ be a positive integer relatively prime to $d$. 
Then, for each positive integer $k\geq 1$, we have the following equality  
$$\# f_{d_1,\cdots, d_n}(k)=\# {\rm Fix}(\sigma^a \circ {\rm Frob}^k|Y({\bar{\bf F}_q})).
\eqno(2.4)$$ 
\end{lemma} 

This lemma was proved in the case $a=1$ in \cite{W3}. It together with the 
general $\ell$-adic fixed point theorem gives 
$$\# f_{d_1,\cdots, d_n}(k)=\sum_{j\geq 0}(-1)^j {\rm Tr}(\sigma \circ {\rm Frob}^k
| H_c^j(Y \otimes {\bar{\rm F}}_q, {\bf Q}_{\ell})),$$
where $\ell$ is a prime number different from $p$ and $H_c^j$ denotes the 
$\ell$-adic cohomology with compact support. This formula is likely    
explicitly stated somewhere in SGA. We have not found it. The quasi-projective case 
is explained in \cite{FW1}. The general finite type case follows by excision.   

Since $\sigma$ and ${\rm Frob}$ commute,  
$\sigma^d=1$, we can decompose the cohomology space into the 
eigenspaces of $\sigma$. The eigenvalues of $\sigma$ are $d$-th roots of unity. 
The eigenvalues of ${\rm Frob}$ are algebraic integers (in fact, Weil $q$-integers  
by Deligne's theorem). 
It follows that there are finitely many $d$-th roots of 
unity $\alpha_i$ and finitely many algebraic integers $\lambda_i$ such that for 
all integers $k\geq 1$, we have  
$$\# f_{d_1,\cdots, d_n}(k) =\sum_i \pm \alpha_i \lambda_i^k.$$
We collect similar terms in terms of 
$\lambda_i$ and rewrite the above expression as  
$$\# f_{d_1,\cdots, d_n}(k) =\sum_j A_j \lambda_j^k,$$
where the $\lambda_j$'s are distinct and $A_j \in {\bf Z}[\zeta_d]$. 
Replacing $\sigma$ by $\sigma^a$ with $(a,d)=1$ and using Lemma 2.3, 
we deduce that for all $\tau\in {\rm Gal}({\bf Q}(\zeta_d)/{\bf Q})$, 
$$\# f_{d_1,\cdots, d_n}(k) =\sum_j \tau(A_j) \lambda_j^k.$$
This sequence of expression is unique since the $\lambda_j$'s are 
distinct. It follows that $\tau(A_j)=A_j$ for all $j$ and all $\tau$. 
Thus, $A_j\in {\bf Z}$ and 
$$Z_{d_1,\cdots, d_n}(f,T) =\prod_j (1-\lambda_j T)^{A_j}$$ 
is indeed a rational function. Theorem 2.2 is proved.

\section{A graph theoretic generalization}

In this section, we give Lenstra's generalization of the partial zeta function 
and its rationality in a graph theory setup. 
Let $G = (V, E)$ be a finite directed graph, where $V$ is the set of 
vertices of $G$ and $E$ is the set of directed edges of $G$. For each edge $e\in E$, let 
$s(e)$ (resp. $t(e)$) denote the starting (resp. the terminal) vertex of the edge $e$. 
Suppose that for each $v\in V$, we are given a scheme $X_v$ of finite type over $\F$. 
Suppose that for each edge $e\in E$, we are given a morphism $f_e: X_{s(e)} \rightarrow X_{t(e)}$ 
of finite type over $\F$. Let $d_v$ ($v\in V$) be positive integers. 
For each positive integer $k$, we define 
$$N(k) =\#\{ x=(x_v)_{v\in V} \in \prod_{v\in V} X({\bf F}_{q^{d_vk}}) | 
\forall e\in E, f_e(x_{s(e)})=x_{t(e)}\}.$$  
Define the graph zeta function to be 
$$Z_{d_1,\cdots, d_n}(G, X, T) = \exp \left(\sum_{k=1}^{\infty} 
{N(k) \over k}T^k \right) \in 1+T{\bf Q}[[T]].$$
One can ask if this power series is a rational function in $T$.  

\begin{theorem}[Lenstra] For any graph $G$, any schemes $X_v$ and any morphisms 
$f_e$ as above, the graph zeta function $Z_{d_1,\cdots, d_n}(G,X,T)$ is   
a rational function in $T$, whose reciprocal zeros and reciprocal poles 
are Weil $q$-integers.  
\end{theorem}

To prove this theorem, it suffices to reduce the above graph zeta function to 
the case of partial zeta functions. For this purpose, let $X$ be the fibred product of the 
schemes $X_v$ ($v\in V$) over all morphisms $f_e$ ($e\in E$). That is, 
$$ X =\{ x\in \prod_{v\in V} X_v | \forall e\in E, f_e(x_{s(e)}) =  x_{t(e)}\}.$$
The scheme $X$ is a closed subscheme of the Cartesian product $\prod_{v\in V} X_v$. 
For each $v\in V$, let $f_v$ be the composed map 
$$f_v:  X \hookrightarrow \prod_{v\in V} X_v \rightarrow X_v,$$
where the last map is the projection to $X_v$. With these definitions, it is clear that 
the graph zeta function $Z_{d_1,\cdots, d_n}(G, X, T)$ is simply the partial 
zeta function $Z_{d_1,\cdots, d_n}(f, T)$ attached to the morphisms $f_v: X\rightarrow X_v$. 
The theorem is proved. 

It may be of interest to explore possible graph theoretic applications of this 
zeta function.

\section{Artin-Schreier hypersurfaces}

To give an example, we consider the case of Artin-Schreier hypersurfaces.  Let 
\[ f(x_1,\dots,x_n,y_1,\dots,y_{n'}) \in {\bf F}_q[x_1,\dots,x_n,y_1,\dots,y_{n'}], \] 
where $n,n' \geq 1$.  For each $d \geq 1$, let 
$$ 
N_d(f) = \#\{(x_0,\dots,x_n,y_1,\dots,y_{n'}):x_0^p-x_0=f(x_1,\dots,x_n,y_1,\dots,y_{n'})\},$$
where $ x_i \in {\bf F}_{q^d}$ ($0\leq i\leq n$) and $ y_j \in {\bf F}_q$ ($1\leq j \leq n'$).
Heuristically (for suitable $f$), we expect 
\[ N_d(f)=q^{dn+n'}+O(q^{(dn+n')/2}) \]
where the constant depends on $p$, $f$, and $d$. Deligne's estimate \cite{De} on 
exponential sums implies the following result. 

\begin{theorem}[Deligne] Given $f$ as above, we 
write $f=f_r+f_{r-1}+\dots+f_0$, where $f_i$ is homogeneous of degree $i$.  
Assume that the 
leading form $f_r$ defines a smooth projective hypersurface 
in ${\bf P}^{n+n'-1}_{{\bf F}_q}$, 
and assume that $p \not| r$.  Then for $d=1$, we have the following inequality 
\[ |N_1(f)-q^{n+n'}| \leq (p-1)(r-1)^{n+n'} q^{(n+n')/2}. \]
\end{theorem}

What can be said about $d>1$? To answer this question, we introduce the 
following terminology.  

\begin{definition}
Let $d$ be a positive integer and let $f$ be a polynomial as above. 
We define the \emph{$d$th fibred sum of $f$} to be the following 
new polynomial 
\[ \textstyle{\bigoplus}_y^d f=f(x_{11},\dots,x_{1n},y_1,\dots,y_{n'})+\dots+f(x_{d1},\dots,x_{dn},y_1,\dots,y_{n'}). \]
\end{definition}

The following estimate on $N_d(f)$ is proved in \cite{FW2}.  

\begin{theorem}[Fu-Wan] Given $f$ as above, we 
write $f=f_r+f_{r-1}+\dots+f_0$, where $f_i$ is homogeneous of degree $i$. 
Assume that $\bigoplus_y^d f_r$ is smooth 
in ${\bf P}^{dn+n'-1}_{{\bf F}_q}$ and assume that $p \not| r$.  Then, 
we have the following inequality 
\[ |N_d(f)-q^{dn+n'}| \leq (p-1)(r-1)^{dn+n'}q^{(dn+n')/2}. \]
\end{theorem}

\begin{example}
Consider the case that 
\[ f(x,y)=f_{1,r}(x_1,\dots,x_n)+f_{2,r}(y_1,\dots,y_{n'})+f_{\leq r-1}(x,y), \]
where $f_{1,r}$ is smooth in ${\bf P}^{n-1}_{{\bf F}_q}$, $f_{2,r}$ is smooth in 
${\bf P}^{n'-1}_{{\bf F}_q}$ and $f_{\leq r-1}$ is a polynomial of degree at most 
$r-1$. It is then straightforward to check that 
$\bigoplus_y^d f_r$ is smooth in ${\bf P}^{dn+n'-1}_{{\bf F}_q}$ 
if and only if $d$ is not divisible by $p$.
Since the condition that the fibred sum be smooth is Zariski open, 
there exist many more examples of such $f$ to which 
the theorem applies if $d$ is not divisible by $p$.
\end{example}

It would be interesting to prove similar results for the Kummer hypersurface 
$x_0^D=f(x_1,\dots,x_n,y_1,\dots,y_{n'})$; see Katz \cite{Ka2} for some related 
weaker results in this direction.

\bigskip

\noindent
{\sl Institute of Mathematics, Chinese Academy of Sciences, Beijing, P.R. China}

\bigskip

\noindent
{\sl Department of Mathematics, University of California, Irvine, CA 92697, USA}

\noindent
{\sl Email: dwan@math.uci.edu}

\end{document}